\documentclass[12pt]{article}
\usepackage{amssymb,amsmath,amsfonts}
\begin{document}

\newtheorem{teo}{Theorem}
\newtheorem{cor}{Corollary}
\newtheorem{lm}{Lemma}
\newtheorem{rem}{Remark}
\newtheorem{pr}{Proposition}
\renewcommand{\theteo}{\arabic{teo}.}
\renewcommand{\thecor}{\arabic{cor}.}
\renewcommand{\thelm}{\arabic{lm}.}
\renewcommand{\therem}{\arabic{rem}.}
\renewcommand{\thepr}{\arabic{pr}.}
\newcommand{\indlim}{\operatornamewithlimits{ind\ lim}}
\renewcommand{\refname}{\begin{center} \mdseries {\sf \large References}
\end{center}}
\begin{center}{\bf ON THE CAUCHY PROBLEM FOR DIFFERENTIAL EQUATIONS IN
A BANACH SPACE OVER THE FIELD OF $p$-ADIC NUMBERS. II.}
\footnote{\noindent Supported by CRDF (Project UM 1-2421-KV-02)}
\end{center}
\vspace*{3mm}
\begin{center} {\textmd \sf MYROSLAV L. GORBACHUK and VALENTYNA I. GORBACHUK}
\end{center}
\vspace*{3mm}

{\bf 1.} \ Let $\mathfrak B$ be a Banach space with norm $\|\cdot\|$ over
the completion $\Omega = {\Omega}_p$ of an algebraic closure of the field
$Q_p$ of $p$-adic numbers ($p$ is prime) (for details we refer to [1 - 3]),
and let $A$ be a closed linear operator on $\mathfrak B$.

For a number $\alpha > 0$, we put
$$
E_{\alpha}(A) = \bigl\{x \in
\bigcap\limits_{n \in {\mathbb N}_0 = \{0, 1, 2, \dots\}}{\mathcal D}(A^n) \Big|
\exists c = c(x) > 0 \quad \forall k \in {\mathbb N}_0 \quad
\|A^k x\| \le c{\alpha}^k\bigr\}
$$
(${\mathcal D}(A)$ is the domain of $A$).

The linear set $E_{\alpha}(A)$ is a Banach space with respect to the norm
$$
\|x\|_{\alpha} = \sup\limits_{n \in {\mathbb N}_0}\frac{\|A^nx\|}{{\alpha}^n}.
$$
Denote by $E(A)$ the space of entire vectors of exponential type for the
operator $A$:
$$
E(A) = \indlim \limits_{\alpha \to \infty} E_{\alpha}(A).
$$
So, as a set,
$$
E(A) = \bigcup\limits_{\alpha > 0} E_{\alpha}(A).
$$

By the type $\sigma (x; A)$ of a vector $x \in E(A)$ we mean the number
$$
\sigma (x; A) = \inf \{\alpha > 0: x \in E_{\alpha}(A)\} =
\varlimsup_{n \to \infty} \|A^nx\|^{\frac{1}{n}}.     \eqno (1)
$$

In the case, where ${\mathcal D}(A) = \mathfrak B$, i.e.,
the operator $A$ is bounded, $E(A) = {\mathfrak B}$, and
$$
\forall x \in {\mathfrak B} \quad \sigma (x; A) \le \|A\|.
$$

{\bf 2.} \ In what follows we shall deal with power series of the form
$$
y(z) = \sum\limits_{n = 0}^{\infty} c_nz^n, \quad c_n \in {\mathfrak B},
\quad z \in \Omega.  \eqno (2)
$$
For such a series the convergence radius is determined by the formula
$$
r = r(y) = \frac{1}{\varlimsup_{n \to \infty}\sqrt[n]{\|c_n\|}}. \eqno (3)
$$
If $r(y) > 0$, then series (2) gives a vector-valued function  $y(z)$
with values in $\mathfrak B$ (a $\mathfrak B$-valued function) in the open disk
$U_r^-(0) = \{z \in \Omega: |z|_p < r\}$ ($|\cdot|_p$ is
the $p$-adic valuation on $\Omega$).

For a number $r > 0$, we denote by ${\mathfrak A}_r(\mathfrak B)$ the set of all
$\mathfrak B$-valued functions $y(z)$ which satisfy the following conditions:

(i) \ $y(z)$ is of the form (2) with $r(y) \ge r$;

(ii) \ $\lim\limits_{i \to \infty}\|c_i\|r^i = 0$.

The linear set ${\mathfrak A}_r(\mathfrak B)$ is a Banach space with respect to
the norm
$$
\|y\|_r = \sup\limits_{n \in {\mathbb N}_0}\|c_n\|r^n.
$$
Moreover, if $0 < r_1 < r$, then the embedding ${\mathfrak A}_r(\mathfrak B) \hookrightarrow
{\mathfrak A}_{r_1}(\mathfrak B)$ induced by restriction of the domain
of a vector-valued function is continuous.

We put
$$
{\mathfrak A}_{loc}(\mathfrak B) = \indlim \limits_{r \to 0}{\mathfrak A}_r(\mathfrak B).
$$
The space ${\mathfrak A}_{loc}(\mathfrak B)$ is called the space of locally analytic
at zero $\mathfrak B$-valued functions.

It follows from (3) that for the convergence radius of $i$-order derivative
$$
y^{(i)}(z) = \sum\limits_{n = 0}^{\infty} (n + 1) \dots (n + i) c_{n + i}z^n,
\quad i \in {\mathbb N},
$$
of a vector-valued function $y(z)$ from ${\mathfrak A}_{loc}(\mathfrak B)$, valid is
the inequality
$$
r(y^{(i)}) \ge r(y).
$$

It is not also hard to check that if $z \to 0$ in $\Omega$, then
$$
y(z) \to y(0) = c_0, \quad \frac{y^{(i)}(z) - y^{(i)}(0)}{z}
\to y^{(i + 1)}(0) = c_{i + 1}(i + 1)!
$$
in the topology of the space $\mathfrak B$.

{\bf 3.} \ Let $x \in E(A), \ \sigma (x; A) = \sigma$. For a fixed natural
$m$ we consider the Mittag-Leffler $\mathfrak B$-valued functions
$$
F_k(z; A)x = \sum\limits_{n = 0}^{\infty} \frac{z^{mn + k}A^nx}
{(mn + k)!}, \quad z \in \Omega, \quad k = 0, 1, \dots, m - 1.  \eqno (4)
$$
\begin{pr} \
The convergence radius of series $(4)$ does not depend on $k$, and it
is determined by the formula
$$
r(F_k(\cdot; A)x) = r = {\sigma}^{-\frac{1}{m}}p^{-\frac{1}{p - 1}}. \eqno (5)
$$
\end{pr}
{\it Proof}. \ According to (3),
$$
r^{-1}(F_k(\cdot; A)x) =
\varlimsup_{n \to \infty}\sqrt[mn + k]{\frac{\|A^nx\|}{|(mn + k)!|_p}}.
$$
By (1), for any $\varepsilon > 0$ and sufficiently large $n \in {\mathbb N}_0$,
$$
\|A^nx\| \le (\sigma + \varepsilon)^n, \eqno (6)
$$
and there exists a subsequence $n_i \to \infty \ (i \to \infty)$ such that
$$
\lim\limits_{i \to \infty}\frac{\|A^{n_i}x\|}{(\sigma - \varepsilon)^{n_i}} =
\infty. \eqno (7)
$$
Using (6) and the estimate
$$
\frac{1}{np}p^{\frac{n}{p - 1}} \le \frac{1}{|n!|_p} \le p^{\frac{n - 1}{p - 1}}
\eqno (8)
$$
valid for large $n \in \mathbb N$ (see [4]), we obtain
$$
\varlimsup_{n \to \infty}\sqrt[mn + k]{\frac{\|A^nx\|}{|(mn + k)!|_p}}
\le (\sigma + \varepsilon)^{\frac{n}{mn + k}}
p^{\frac{mn + k - 1}{(mn + k)(p - 1)}}.
$$
Since
$$
\lim\limits_{n \to \infty}\,(\sigma + \varepsilon)^{\frac{n}{mn + k}}
p^{\frac{mn + k - 1}{(mn + k)(p - 1)}} = (\sigma + \varepsilon)^{\frac{1}{m}}
p^{\frac{1}{(p - 1)}},
$$
the inequality
$$
\varlimsup_{n \to \infty}\sqrt[mn + k]{\frac{\|A^nx\|}{|(mn + k)!|_p}}
\le (\sigma + \varepsilon)^{\frac{1}{m}}
p^{\frac{1}{(p - 1)}}
$$
holds. Taking into account that $\varepsilon > 0$ is arbitrary, we conclude
that
$$
\varlimsup_{n \to \infty}\sqrt[mn + k]{\frac{\|A^nx\|}{|(mn + k)!|_p}}
\le {\sigma}^{\frac{1}{m}}p^{\frac{1}{(p - 1)}},
$$
whence
$$
r(F_k(\cdot; A)x) \ge {\sigma}^{-\frac{1}{m}}p^{-\frac{1}{(p - 1)}}.
$$

If now $\{n_i\}_{i = 1}^{\infty}$ is a subsequence satisfying (7), then because
of (8),
$$
\sqrt[mn_i + k]{\frac{\|A^{n_i}x\|}{|(mn + k)!|_p}} \ge
(\sigma - \varepsilon)^{\frac{n_i}{mn_i + k}}p^{\frac{1}{(p - 1)}}
[p(mn_i + k)]^{-\frac{1}{(mn_i + k)}}.
$$
Since
$$
\lim\limits_{n \to \infty}\,(\sigma - \varepsilon)^{\frac{n_i}{mn_i + k}}
[p(mn_i + k)]^{-\frac{1}{(mn_i + k)}} = (\sigma - \varepsilon)^{\frac{1}{m}},
$$
we have
$$
\varlimsup_{n \to \infty}\sqrt[mn + k]{\frac{\|A^nx\|}{|(mn + k)!|_p}}
\ge (\sigma - \varepsilon)^{\frac{1}{m}}p^{\frac{1}{(p - 1)}},
$$
whence
$$
r(F_k(\cdot; A)x) \le {\sigma}^{-\frac{1}{m}}p^{-\frac{1}{(p - 1)}}.
$$
Thus, formula (5) is valid. Q.E.D..

One can easily check that for the derivatives of $F_k(z; A)x)$, the
relations
$$
F_k^{(i)}(z; A)x) =
\left\{
\begin{array}{rcl}
F_{k - i}(z; A)x & \mbox{\rm if} & 0 \le i \le k \\
F_0^{(i - k)}(z; A)x & \mbox{\rm if} & i > k \\
\end{array} \right.  \eqno (9)
$$
and
$$
\left\{
\begin{array}{rcl}
F_0^{(ml)}(z; A)x & = & F_0(z; A)A^lx \\
F_0^{(ml + j)}(z; A)x & = & (F_0(z; A)A^lx)^{(j)} = F_{m - j}(z; A)A^{l + 1}x
\ 1 \le j \le m - 1 \\
\end{array} \right.  \eqno (10)
$$
are fulfilled, which imply, with regard to $\sigma(A^jx; A) = \sigma(x; A) \
(j \in \mathbb N)$, that
$$
r(F_k^{(i)}(\cdot; A)x) = {\sigma}^{-\frac{1}{m}}p^{-\frac{1}{(p - 1)}}
\quad i \in {\mathbb N}_0, \ k = 0, 1, \dots, m - 1.
$$

{\bf 4.} \ Let us consider the Cauchy problem
$$
\left\{
\begin{array}{rcl}
y^{(m)}(z) & = & Ay(z) \\
y^{(k)}(0) & = & y_k, \ \ k = 0, 1, \dots, m - 1, \\
\end{array} \right.  \eqno (11)
$$
where $A$ is a closed linear operator on $\mathfrak B$ (the case $m = 1$ was
discussed in [5]). By a solution of this
problem we mean a vector-valued function of the form (2) with values in
${\mathcal D}(A)$ that satisfies (11). The question arises, under what conditions
on the initial data $y_k$ problem (11) has a solution in the space
${\mathfrak A}_{loc}(\mathfrak B)$. The following assertion gives an answer.
\begin{teo} \
Problem $(11)$ is solvable in the class ${\mathfrak A}_{loc}(\mathfrak B)$ if and only
if $y_k \in E(A), \ k = 0, 1, \dots, m - 1$. The solution is represented
in the form
$$
y(z) = \sum\limits_{k=0}^{m-1} F_k(z, A)y_k.  \eqno (12)
$$
Moreover, problem $(11)$ is well-posed in ${\mathfrak A}_{loc}(\mathfrak B)$, that is,
the solution is unique in this space, and the convergence $E(A) \ni y_{i,k} \to y_k \
(i \to \infty, \ k = 0, 1, \dots, m - 1)$ in the $E(A)$-topology implies
the convergence of the sequence of corresponding solutions $y_i(z)$ to
$y(z)$ in ${\mathfrak A}_{loc}(\mathfrak B)$.
\end{teo}
{\it Proof}. \  Assume that problem (11) is solvable in
${\mathfrak A}_{loc}(\mathfrak B)$. This means that there exists a
${\mathcal D}(A)$-valued function $y(z) = \sum\limits_{n = 0}^{\infty}
c_nz^n \in {\mathfrak A}_r(\mathfrak B)$ with some $r > 0 \ (c_n \in {\mathfrak B},
\ z \in \Omega)$, which satisfies (11). It is obvious that $y_0 = y(0) \in
{\mathcal D}(A)$. Further, if we take $|\Delta z|_p < |z|_p$, then
$|z + \Delta z|_p = |z|_p$, the inequality $|z|_p < r(y)$ implies
$|z + \Delta z|_p < r(y)$, and for any $z \in U_r^-(0)$ we have
$$
\frac{y^{(m)}(z + \Delta z) - y^{(m)}(z)}{\Delta z} =
A\frac{y(z + \Delta z) - y(z)}{\Delta z}.
$$
In view of closedness of the operator $A$ and the fact [3]
that
$$
\frac{y^{(k)}(z + \Delta z) - y^{(k)}(z)}{\Delta z} \to y^{(k + 1)}(z),
$$
the values of $y'(z)$ in $U_r^-(0)$ belong to ${\mathcal D}(A)$, and
$$
y^{(m + 1)}(z) = A y'(z) \quad \mbox{\rm when} \ |z|_p < r(y).
$$
Repeating this procedure as much as we need, we arrive at the inclusion that
$$
y^{(n)}(z) \in {\mathcal D}(A) \ \ \mbox{\rm if} \ \ |z|_p < r(y) \le r, \quad
n \in {\mathbb N}_0,
$$
and at the equality
$$
y^{(mn + k)}(z) = A^n y^{(k)}(z), \quad |z|_p < r(y) \le r,
\quad k = 0, 1, \dots, m - 1,
$$
whence
$$
A^ny_k = y^{(mn + k)}(0) = (mn + k)!c_{mn + k} \eqno (13)
$$
By property (ii) of a $\mathfrak B$-valued function $y \in {\mathfrak A}_r(\mathfrak B)$
and the estimate (8),
$$
\|A^ny_k\| = |(mn + k)!|_p \|c_{mn + k}\| \le
\frac{(mn + k)p\|c_{mn + k}\|p^{-\frac{mn + k}{p - 1}}}{r^{mn + k}}r^{mn + k}
\le c{\alpha}^n,
$$
where
$$
c = p\left(\frac{(1 + \varepsilon)p^{-\frac{1}{p - 1}}}{r}\right)^k
\sup\limits_{i}\,\|c_i\|r^i, \quad \alpha =
\left(\frac{(1 + \varepsilon)p^{-\frac{1}{p - 1}}}{r}\right)^m
$$
($\varepsilon > 0$ is arbitrary). So, $y_k \in E(A), \ k = 0, 1, \dots, m - 1$.

Conversely, let $y_k \in E(A), \ k = 0, 1, \dots, m - 1$. Taking into account
relations (9) and (10), one can verify that the $\mathfrak B$-valued functions
$F_k(z; A)x, \ x \in E(A),$ satisfy the equation from (11) in $U_r^-(0)$
with $r = {\sigma(x; A)}^{-\frac{1}{m}}p^{-\frac{1}{p - 1}}$ and the initial
data
$$
F_k^{(i)}(0; A)x = {\delta}_{ik}x, \quad  i, k = 0, 1, \dots, m - 1.
$$
Therefore the ${\mathcal D}(A)$-valued function (12) is a solution of (11) from
the space ${\mathfrak A}_{loc}(\mathfrak B)$.
It follows from (13) that this solution is unique in the mentioned class.

To prove the well-posedness of problem (11) in ${\mathfrak A}_{loc}(\mathfrak B)$,
assume that a sequence $y_{i,k} \in E(A)$ converges in $E(A)$ to $y_k$ as
$i \to \infty, \ k = 0, 1, \dots, m - 1$. This means that there exists
$\alpha > 0$ such that $y_{i,k}, y_k \in E_{\alpha}(A)$, and
$$
\|y_{i,k} - y_k\|_{\alpha} \to 0 \ \ \mbox{\rm as} \ i \to \infty.
$$
Then the corresponding solutions $y_i(z)$ and $y(z)$ of problem (11)
belong to ${\mathfrak A}_r(\mathfrak B)$ where
$r < {\alpha}^{-\frac{1}{m}}p^{-\frac{1}{p - 1}}$, and
$$
\|y_i(\cdot) - y(\cdot)\|_r \le
\sum\limits_{k = 0}^{m - 1}\|F_k(\cdot; A)y_k\|_r =
\sum\limits_{k = 0}^{m - 1} \sup\limits_{n \in {\mathbb N}_0}
\left\|\frac{A^n(y_k - y_{i,k})}{(mn + k)!}
\right\| r^{mn + k}
\le
$$
$$
\sum\limits_{k = 0}^{m - 1} \sup\limits_{n \in
{\mathbb N}_0} \frac{{\alpha}^n\|y_k - y_{i,k}\|_{\alpha}
{\alpha}^{-\frac{mn + k}
{m}}p^{-\frac{mn + k}{p - 1}}}{|(mn + k)!|_p} \le
$$
$$
\sum\limits_{k = 0}^{m - 1} \sup\limits_{n \in {\mathbb N}_0}
{\alpha}^{-\frac{k}{m}}p^{-\frac{mn + k}{p - 1}}
p^{\frac{mn + k - 1}{p - 1}}\|y_{i,k} - y_k\|_{\alpha} =
$$
$$
\sum\limits_{k = 0}^{m - 1} {\alpha}^{-\frac{k}{m}}
p^{-\frac{1}{p - 1}}\|y_{i,k} - y_k\|_{\alpha} = p^{-\frac{1}{p - 1}}
\sum\limits_{k = 0}^{m - 1} {\alpha}^{-\frac{k}{m}}
\|y_{i,k} - y_k\|_{\alpha}.
$$
Since $\|y_{i,k} - y_k\|_{\alpha} \to 0$
as $i \to \infty \ (k = 0, 1, \dots, m - 1)$, we have
$$
\|y_i(\cdot) - y(\cdot)\|_r \to 0 \ \ (i \to \infty),
$$
which completes the proof.

It follows from the above proof, that in order that problem (11) be solvable
in the class of entire $\mathfrak B$-valued functions, it is necessary and
sufficient that
$$
y_k \in \bigcap\limits_{\alpha > 0} E_{\alpha}(A).
$$
\begin{cor} \
If the operator $A$ is bounded, then problem $(11)$ is well-posed in
${\mathfrak A}_r(\mathfrak B)$ for any $y_k \in \mathfrak B$.
\end{cor}
\begin{rem} \
As is shown in $[6]$, in the case where $\mathfrak B$ is a Banach space over the field
$\mathbb C$ of complex numbers, the Cauchy problem $(11)$ is well-posed
in the class ${\mathfrak A}_r(\mathfrak B)$ if and only if
$$
\forall \alpha > 0 \quad \exists c = c(\alpha) > 0 \quad \|A^ny_k\| \le c
{\alpha}^nn^{mn}, \quad k = 0, 1, \dots, m - 1.
$$
In order that the solution of $(11)$ be an entire $\mathfrak B$-valued function of
exponential type, it is necessary and sufficient that $y_k \in E(A), \
0, 1, \dots, m - 1$.
\end{rem}

{\bf 5.} \ In this section we show how the result of Theorem 1 can be applied
to partial differential equations (see also [7]).

Let ${\mathcal A}_{\rho}$ be the space of $\Omega$-valued functions $f(x)$ analytic
on the $n$-dimensional disk
$$
U_{\rho}^+(0) = \left\{x = (x_1, \dots x_n) \in {\Omega}^n:
|x|_p = \left(\sum\limits_{i = 1}^{n}|x_i|_p^2\right)^{1/2}
\le \rho\right\}.
$$
This means that
$$
f(x) = \sum\limits_{\alpha} f_{\alpha}x^{\alpha}, \quad f_{\alpha} \in \Omega,
\quad \lim\limits_{|\alpha| \to \infty}|f_{\alpha}|_p{\rho}^{|\alpha|} = 0,
$$
where $\alpha = ({\alpha}_1, \dots
{\alpha}_n), \ {\alpha}_i \in {\mathbb N}_0, \ |\alpha| = {\alpha}_1 + \dots + {\alpha}_n$.

The space ${\mathcal A}_{\rho}$ is a $p$-adic Banach space with respect
to the norm
$$
\|f\|_{\rho} = \sup\limits_{\alpha}|f_{\alpha}|_p {\rho}^{|\alpha|}.
$$
It is clear that the differential operators
$$
f \mapsto \frac{\partial f}{\partial x_j} = \sum\limits_{\alpha} {\alpha}_j f_{\alpha}
x_1^{{\alpha}_1} \dots x_j^{{\alpha}_j - 1} \dots x_n^{{\alpha}_n}, \ j =
1, \dots, n,
$$
are bounded in ${\mathcal A}_{\rho}$, and
$$
\left\|\frac{\partial f}{\partial x_j}\right\|_{\rho} =
\sup\limits_{\alpha}|f_{\alpha}{\alpha}_j|_p{\rho}^{|\alpha| - 1} \le
\frac{1}{\rho}\sup\limits_{\alpha}|f_{\alpha}|_p{\rho}^{|\alpha|} =
\frac{1}{\rho}\|f\|_{\rho}.
$$
Since for $f(x) = x_j, \
\left\|\frac{\partial f}{\partial x_j}\right\|_{\rho} = \frac{1}{\rho}$, the
norm of the operator $\frac{\partial}{\partial x_j}$ is equal to
$\frac{1}{\rho}$.

The multiplication operator
$$
G: f \mapsto fg, \ f \in {\mathcal A}_{\rho}, \quad g \in {\mathcal A}_{\rho},
$$
is bounded in ${\mathcal A}_{\rho}$, too, and
$$
\|G\| = \|g\|_{\rho}.
$$
Indeed, let $f(x) = \sum\limits_{\alpha} f_{\alpha}x^{\alpha}, \
g(x) = \sum\limits_{\alpha} g_{\alpha}x^{\alpha}, \ \alpha = ({\alpha}_1,
\dots, {\alpha}_n)$. Then
$$
f(x)g(x) = \sum\limits_{\alpha} c_{\alpha}x^{\alpha},
$$
where
$$
c_{\alpha} = \sum\limits_{0 \le i \le \alpha} f_ig_{\alpha - i}
= \sum\limits_{i_1 = 0}^{{\alpha}_1} \dots \sum\limits_{i_n = 0}^{{\alpha}_n}
f_{i_1, \dots, i_n}g_{{\alpha}_1 - i_1, \dots, {\alpha}_n - i_n} \quad
(i = (i_1, \dots, i_n)).
$$
So,
$$
\|fg\|_{\rho} = \sup\limits_{\alpha}\max_{0 \le i \le \alpha} |f_i|_p
|g_{\alpha - i}|_p {\rho}^{|i|}{\rho}^{|\alpha| - |i|} \le \|f\|_{\rho}
\|g\|_{\rho}.
$$
As $\|fg\|_{\rho} = \|g\|_{\rho}$ for $f \equiv 1$, we have $\|G\| =
\|g\|_{\rho}$.

Let us consider now the Cauchy problem
$$
\left\{
\begin{array}{rcl}
{\displaystyle \frac{{\partial}^m u(t, x)}{\partial t^m}} & = &
\sum\limits_{|\beta| = 0}^na_{\beta}(x)D^{\beta}u(t, x) \\
u^{(k)}(0, x) & = & {\varphi}_k (x), \ \ k = 0, 1, \dots, m - 1, \\
\end{array}
\right. \eqno (14)
$$
where
$$
a_{\beta}(x) \in {\mathcal A}_{\rho}, \quad {\varphi}_k (x) \in {\mathcal A}_{\rho},
\quad D^{\beta} = {\displaystyle \frac{{\partial}^{|\beta|}}
{\partial x_1^{{\beta}_1} \dots \partial x_n^{{\beta}_n}}}.
$$

If we put ${\mathfrak B} = {\mathcal A}_{\rho}$ and define the operator $A$
as
$$
f \mapsto Af = \sum\limits_{|\beta| = 0}^na_{\beta}D^{\beta}f,
$$
then problem (14) can be written in the form (11). Furthermore, the relations
$$
\bigl\|\sum\limits_{|\beta| = 0}^na_{\beta}D^{\beta}f\bigr\|_{\rho} \le
\max\limits_{\beta}\|a_{\beta}D^{\beta}f\|_{\rho} \le
\max\limits_{\beta}\|a_{\beta}\|_{\rho}\|D^{\beta}f\|_{\rho} \le
\max\limits_{\beta}\{{\rho}^{-|\beta|}\|a_{\beta}\|_{\rho}\}\|f\|_{\rho}
$$
show that the operator $A$ is bounded in ${\mathcal A}_{\rho}$, and
$$
\|A\| \le \max\limits_{\beta}\{{\rho}^{-|\beta|}\|a_{\beta}\|_{\rho}\}.
$$
It follows from Corollary 1 that problem (14) is well-posed in
${\mathfrak A}_{loc}({\mathcal A}_{\rho})$ in the disk
$t \in \Omega: |t|_p <
p^{-\frac{1}{p - 1}}(\max\limits_{\beta}\,{\rho}^{-|\beta|}
\|a_{\beta}\|_{\rho})^{-\frac{1}{m}}$.
\newpage
%\bigskip
%\begin{center}{REFERENCES} \end{center}
%\smallskip

\vspace*{3mm}
Institute of Mathematics \\
National Academy of Sciences of Ukraine \\
3 Tereshchenkivs'ka \\
Kyiv 01601, Ukraine \\
E-mail: imath@horbach.kiev.ua

\end{document}